\documentclass[leqno,12pt]{amsart}
\usepackage{latexsym}
\usepackage{mathrsfs}
\usepackage{amsmath}
\usepackage{amssymb}
\usepackage[bookmarks]{hyperref}
\usepackage{pstricks,pst-plot}

\font\tenscr=rsfs10 
\font\sevenscr=rsfs7 
\font\fivescr=rsfs5 
\skewchar\tenscr='177 \skewchar\sevenscr='177 \skewchar\fivescr='177
\newfam\scrfam \textfont\scrfam=\tenscr \scriptfont\scrfam=\sevenscr
\scriptscriptfont\scrfam=\fivescr
\def\scr{\fam\scrfam}

\newtheorem{theorem}{Theorem}[section]
\newtheorem{lemma}[theorem]{Lemma}
\newtheorem{corollary}[theorem]{Corollary}
\newtheorem{proposition}[theorem]{Proposition}

\theoremstyle{definition}
\newtheorem{remark}[theorem]{Remark}
\newtheorem{definition}[theorem]{Definition}

\allowdisplaybreaks[1]

\newcommand{\C}{\mathbb{C}}
\def\C{{\mathbf {C}\/}}

\newcommand{\R}{\mathbf{R}}
\def\R{{\mathbf {R}\/}}

\newcommand{\row}[2]{#1_1,\ldots,#1_#2}
\newcommand{\range}[1]{1,\ldots, #1}

\def\Cn{{\C^N}}
\def\rhull {h_r(X)}

\def\row#1#2{#1_1,\ldots,#1_#2}

\newcommand{\what}{\widehat}

\def\wX{\widehat X}

\def\cma{{\frak M}_A}

\newcommand{\bthm}{\begin{theorem}}
\newcommand{\ethm}{\end{theorem}}
\newcommand{\blem}{\begin{lemma}}
\newcommand{\elem}{\end{lemma}}
\newcommand{\bcor}{\begin{corollary}}
\newcommand{\ecor}{\end{corollary}}
\newcommand{\bprop}{\begin{proposition}}
\newcommand{\eprop}{\end{proposition}}
\newcommand{\bdefn}{\begin{definition}}
\newcommand{\edefn}{\end{definition}}
\newcommand{\bpf}{\begin{proof}}
\newcommand{\epf}{\end{proof}}

\newcommand{\bm}{\bibitem}

\newcommand{\bi}{\begin{itemize}}
\newcommand{\ei}{\end{itemize}}

\newcommand{\bc}{\begin{cases}}
\newcommand{\ec}{\end{cases}}

\newcommand{\ba}{\begin{array}}
\newcommand{\ea}{\end{array}}

\newcommand{\be}{\begin{equation}}
\newcommand{\ee}{\end{equation}}

\newcommand{\bea}{\begin{eqnarray}}
\newcommand{\eea}{\end{eqnarray}}

\newcommand{\beaa}{\begin{eqnarray*}}
\newcommand{\eeaa}{\end{eqnarray*}}

\newcommand{\beastar}{\begin{eqnarray*}}
\newcommand{\eeastar}{\end{eqnarray*}}

\font\tenscr=rsfs10 
\font\sevenscr=rsfs7 
\font\fivescr=rsfs5 
\skewchar\tenscr='177 \skewchar\sevenscr='177 \skewchar\fivescr='177
\newfam\scrfam \textfont\scrfam=\tenscr \scriptfont\scrfam=\sevenscr
\scriptscriptfont\scrfam=\fivescr
\def\scr{\fam\scrfam}

\def\scr{\fam\scrfam}

\def\Bone#1{{\scr B}(#1)}
\def\Bzero#1{{\scr B}_0(#1)}

\def\forallBzero#1{\hbox{for all } f\in {{\scr B}_0(#1)}}
\def\forallBone#1{\hbox{for all } f\in {{\scr B}(#1)}}

\def\textfrac#1#2{{\textstyle \frac{#1}{#2}}}

\def\pB{\partial B}
\def\ob{\overline B}
\def\tB{\tilde B}
\def\od{\overline D}
\def\disc#1#2#3{D(#1_{#3},#2_{#3})} 
\def\odisc#1#2#3{{\overline D}(#1_{#3},#2_{#3})}
\def\pfz{{{\partial f}/{\partial z_l}}}
\def\oh{{\overline H}}
\def\endhat{\widehat{\phantom |}}

\def\varep{\varepsilon}

\def\pji{p_j^{-1}}
\def\pjid#1{p_j^{-1}(#1 D)}

\def\inv{\mathop{\rm inv}\nolimits}

\begin{document}

\title[Uniform algebras with dense invertible group]{Gleason parts and point derivations\\ for uniform algebras with\\ dense invertible group}
\author{Alexander J. Izzo}
\address{Department of Mathematics and Statistics, Bowling Green State University, Bowling Green, OH 43403}
\email{aizzo@bgsu.edu}



\begin{abstract}
It is shown$\vphantom{\widehat{\widehat{\widehat{\widehat{\widehat{\widehat{\wX}}}}}}}$ that there exists a compact set $X$ in $\C^N$ ($N\geq 2$) such that $\wX\setminus X$ is nonempty and the uniform algebra $P(X)$ has a dense set of invertible elements, a large Gleason part, and an abundance of nonzero bounded point derivations.  The existence of a Swiss cheese $X$ such that $R(X)$ has a Gleason part of full planar measure and a nonzero bounded point derivation at almost every point is established.  An analogous result in $\C^N$ is presented.  The analogue for rational hulls of a result of Duval and Levenberg on polynomial hulls containing no analytic discs is established.  The results presented address questions raised by Dales and Feinstein.
 
\end{abstract}
\maketitle

\vskip -2.55 true in
\centerline{\footnotesize\it Dedicated to Andrew Browder} 
\vskip 2.55 truein

\section{Introduction}

It was once conjectured that whenever the polynomially convex hull $\wX$ of a compact set $X$ in $\C^N$ is strictly larger than $X$, the complementary set $\wX\setminus X$ must contain an analytic disc.  This conjecture was disproved by
Gabriel~Stolzenberg \cite{Stol}.  Given Stolzenberg's result, it is natural to ask whether weaker semblances of analyticity, such as nontrivial Gleason parts or nonzero bounded point derivations, must be present in $\wX\setminus X$.  
Garth~Dales and Joel~Feinstein \cite{DalesF} strengthened Stolzenberg's result by constructing a compact set $X$ in $\C^2$ with $\wX\setminus X$ nonempty and such that $P(X)$ has a dense set of invertible elements.  (As noted in \cite{DalesF} the condition that $P(X)$ has a dense set of invertible elements is strictly stronger than the condition that $\wX$ contains no analytic disc.)  In the paper \cite{CGI} of Brian~Cole, Swarup~Ghosh, and Alexander~Izzo, an example is given in $\C^3$ with the additional properties that $P(X)$ has no nontrivial Gleason parts and no nonzero bounded point derivations.  The main purpose of the present paper is to show that in contrast to that example, there are also examples having a large Gleason part and an abundance of nonzero bounded point derivations.  In particular we will establish the following result.  Here $\partial B$ denotes the boundary of the open unit ball $B=\{z: \|z\|<1\}$ in $\C^N$, and following Dales and Feinstein, 
we will say that a uniform algebra $A$ has dense invertibles if the invertible elements of $A$ are dense in $A$.

\begin{theorem}\label{maintheorem1}
For each integer $N\geq 2$, there exists a compact set $X\subset \partial B\subset \Cn$ such that
\item{\rm(i)} $\wX\setminus X$ has positive $2N$-dimensional measure
\item{\rm(ii)} $P(X)$ has dense invertibles
\item{\rm(iii)} there is a set $P\subset \wX\setminus X$ of positive 
$2$-dimensional Hausdorff measure contained in a single Gleason part for $P(X)$ 
\item{\rm(iv)} at each point of $P$ there is a nonzero bounded point derivation for $P(X)$.
\end{theorem}

In the case $N\geq 3$ we can obtain a larger Gleason part and more bounded point derivations, but at the expense of giving up assurance that $\wX\setminus X$ has positive $2N$-dimensional measure.

\begin{theorem}\label{maintheorem2}
For each integer $N\geq 3$, there exists a compact set $X\subset \partial B\subset \Cn$ such that
\item{\rm(i)} $P(X)$ has dense invertibles
\item{\rm(ii)} there is a set $P\subset \wX\setminus X$ of positive 
$2(N-1)$-dimensional Hausdorff measure contained in a single Gleason part for $P(X)$ 
\item{\rm(iii)} at each point of $P$ the space of bounded point derivations for $P(X)$ has dimension at least~$N-1$.
\end{theorem}

One can also consider questions of analyticity for the rationally convex hull $h_r(X)$ of a compact set $X$.  In fact, Theorem~\ref{maintheorem2} will be obtained as an immediate consequence of the next result which concerns rationally convex hulls.

\begin{theorem}\label{maintheorem3}
For each integer $N\geq 2$, there exists a compact set $X\subset \partial B\subset \Cn$ such that 
\item{\rm(i)} $R(X)$ has dense invertibles
\item{\rm(ii)} there is a set $P\subset \rhull\setminus X$ of positive $2N$-dimensional measure contained in a single Gleason part for $R(X)$ 
\item{\rm(iii)} at each point of $P$ the space of bounded point derivations for $R(X)$ has dimension~$N$.
\end{theorem}

The proofs of Theorems~\ref{maintheorem1} and~\ref{maintheorem2} show that the set $X$ in those theorems satisfies $\wX=\rhull$.  We have not stated this in the theorems because, as noted by Dales and Feinstein, this condition actually follows from the condition that $P(X)$ has dense invertibles.

It is surprising that little attention has been given in the literature to the existence of nontrivial Gleason parts in hulls without analytic discs because early efforts to prove the existence of analytic discs focused on the use of Gleason parts.  
John~Wermer \cite{Wermer} proved that every nontrivial Gleason part for a Dirichlet algebra is an analytic disc, and this was extended to uniform algebras with uniqueness of representing measures by Gunter~Lumer \cite{Lumer}.  However, once it was proven that analytic discs  do not always exist, it seems that almost no further work was done on Gleason parts in hulls except for a result of Richard~Basener \cite{Basener} that there exists a compact set in $\C^2$ whose rationally convex hull contains no analytic discs but contains a Gleason part of positive 
4-dimensional measure.  Theorem~\ref{maintheorem3} above strengthens this result.  Our method of proof is quite different from Basener's.

In connection with their result mentioned above,
Dales and Feinstein asked what can be said about the existence of nonzero point derivations \cite[Section~4, Question~1]{DalesF}.  
Although they asked specifically about the algebra $P(Y)$ for the set $Y$ that they constructed, the spirit of their question is what can be said about point derivations on uniform algebras of the form $P(X)$ ($X\subset\C^N$ compact) at points of the set $\wX\setminus X$, when $X$ is such that $\wX\setminus X$ is nonempty while $P(X)$ has dense invertibles.  In fact, the question of existence of point derivations in hulls without analytic discs seems not to have been addressed in the literature at all until very recently.
Taken together, Theorem~\ref{maintheorem1} above and \cite[Theorem~1.4]{CGI} answer this question in the case of bounded point derivations: there exist such uniform algebras both with, and without, bounded point derivations.  

Dales and Feinstein \cite[Section~4, Question~2]{DalesF} also raised the question of whether the set of exponentials $\exp A$ can be dense in a uniform algebra $A$ with proper Shilov boundary (so in particular in a $P(X)$ with $\wX\setminus X\not= \emptyset$).  They also raised the question of whether $\exp A$ can ever be dense in a nontrivial uniform algebra $A$.  In connection with these questions, they remarked that \lq\lq it is conceivable (but unlikely)\rq\rq\ that $\exp A$ is dense in the algebra $A$ with $A=P(Y)$ where $Y$ is the set they constructed.  We will establish the following.

\bthm \label{exp}
The examples in Theorems~\ref{maintheorem1}--\ref{maintheorem3} can be taken to have the additional property that $\exp A$ is {\em not\/} dense in $A$ where $A=P(X)$ in Theorems~\ref{maintheorem1} and~\ref{maintheorem2} and $A=R(X)$ in Theorem~\ref{maintheorem3}.
\ethm

This, of course, leaves unsettled the question of whether there exist nontrivial uniform algebras for which the exponentials {\em are\/} dense.

It is also of interest to consider Gleason parts and point derivations in algebras with dense invertibles without the requirement that the hull be nontrivial.
For uniform algebras on planar sets some results on this have long been known.  For a compact set $X\subset \C$, it follows from Lavrentiev's theorem \cite[Theorem~3.4.15]{Browder} that if $P(X)$ has dense invertibles then $P(X)=C(X)$ and consequently has no nonzero point derivations \cite[pp.~64--65]{Browder} and no nontrivial Gleason parts.  With regard to $R(X)$,
interesting examples concerning point derivations have been given using Swiss cheeses.  The definition of a Swiss cheese we shall use is as follows: A 
{\it Swiss cheese\/} is a compact set $X$ obtained from the closed unit disc 
$\od$ by deleting a sequence of open discs $\{D_j\}_{j=1}^\infty$ with radii 
$\{r_j\}_{j=1}^\infty$ such that the closures of the $D_j$ are disjoint, 
$\sum_{j=1}^\infty r_j<\infty$, and the resulting set $X=\od\backslash(\cup_{j=1}^\infty D_j)$ has no interior.  (Sometimes in the literature a more general definition of Swiss cheese is used and Swiss cheeses as we have defined them are referred to as {\em classical Swiss cheeses\/}.)  It is well known that every Swiss cheese has positive 2-dimensional measure and satisfies $R(X)\not=C(X)$.  Andrew Browder constructed a  Swiss cheese $X$ such that $R(X)$ has a bounded point derivation at almost every point with respect to 2-dimensional measure (see \cite{Wermerderivation}), and Wermer constructed a Swiss cheese $X$ such that $R(X)$ has no nonzero bounded point derivations \cite{Wermerderivation}.    We will prove the following result which deals with both bounded point derivations and Gleason parts.

\begin{theorem}\label{cheese}
There exists a Swiss cheese $X$  such that the set of nonpeak points for $R(X)$ is a single Gleason part of full 2-dimensional measure in $X$, and $R(X)$ has a bounded point derivation at almost every point.
\end{theorem}

The question of whether there exists a Swiss cheese $X$ for which $R(X)$ is antisymmetric was raised by Feinstein (private communication).  One might expect that every Swiss cheese would have this property, but a counterexample was constructed by Lynn Steen \cite{Steen}.  We will show that the Swiss cheese of Theorem~\ref{cheese} 
provides an affirmative answer to Feinstein's question.

We extend Theorem~\ref{cheese} to $N$ dimensions as follows.

\begin{theorem}\label{cheese-dim-n}
There exists a compact rationally convex set $X\subset \C^N$ {\rm(}$N\geq 1${\rm)} of positive $2N$-dimensional measure such that $R(X)$ has dense invertibles, there is a Gleason part for $R(X)$ of full $2N$-dimensional  measure in $X$,  and at almost every point {\rm(}with respect to $2N$-dimensional measure{\rm)}  the space of bounded point derivations for $R(X)$ has dimension $N$.
\end{theorem}

In fact, the proof will show that $X$ can be taken to be a subset of the ball $\ob$ such that the complement of $X$ in $\ob$ has $2N$-dimensional measure less than any prescribed $\varep>0$.
The theorem can be translated into a statement about polynomially convex sets in $\C^{N+1}$ by a standard device recalled in the next section.  This is left to the reader.

The proof of Theorem~\ref{maintheorem1} is based on the approach to constructing hulls without analytic discs due to Julien~Duval and Norman~Levenberg~\cite{DuvalL}.  Their result is as follows except that here we make the additional observation that the polynomial hull and rational hull of the set constructed coincide.

\bthm \label{DL}
If $K$ is a compact polynomially convex subset of $B\subset \C^N$ {\rm(}$N\geq 2${\rm)}, then there is a compact subset $X$ of $\pB$ such that $\wX=\rhull\supset K$ and the set $\wX\setminus K$ contains no analytic disc.
\ethm

Theorem~\ref{maintheorem1} will be proven by using the method of Duval and Levenberg to construct a compact set $X\subset \pB$ whose polynomial hull contains a graph over a Swiss cheese.  The presence of this subset in $\wX$ gives the existence of the desired Gleason part and bounded point derivations.
To prove Theorem~\ref{maintheorem3} a different construction that enables us to put a specified {\it rationally\/} convex set into a rational hull is needed.  This leads us to prove the following analogue of the result of Duval and Levenberg with polynomial convexity replaced by rational convexity.

\bthm \label{rationalDL}
If $K$ is a compact rationally convex subset of $B\subset \C^N$ {\rm(}$N\geq 2${\rm)}, then there is a compact subset $X$ of $\partial B$ such that $\rhull\supset K$ and the set $\rhull\setminus K$ contains no analytic disc.
\ethm

In the next section we recall some definitions and notations already used above.  
Section~3 is devoted to the theorem of Duval and Levenberg and its analogue for rational convexity.  The hulls with dense invertibles and large Gleason parts and nonzero bounded point derivations (Theorems~\ref{maintheorem1}--\ref{maintheorem3}) are constructed in Section~4.  The Swiss cheese of Theorem~\ref{cheese} is constructed in Section~5, and the $N$-dimensional extension given by Theorem~\ref{cheese-dim-n} is presented in Section~6.  The paper concludes with some open questions in Section~7.

It is a pleasure to dedicate this paper to Andrew Browder.  As a graduate student first learning about uniform algebras, I found his book \cite{Browder} immensely helpful and inspiring.  When I met Browder I found him to be very friendly and  generous with his knowledge and ideas.  Over the years I have appreciated his insightful observations regarding my work.


\section{Preliminaries}~\label{prelim}

For
$X$ a compact (Hausdorff) space, we denote by $C(X)$ the algebra of all continuous complex-valued functions on $X$ with the supremum norm
$ \|f\|_{X} = \sup\{ |f(x)| : x \in X \}$.  A \emph{uniform algebra} on $X$ is a closed subalgebra of $C(X)$ that contains the constant functions and separates
the points of $X$.  

For a compact set $X$ in $\C^N$, we denote by 
$P(X)$ the uniform closure on $X$ of the polynomials in the complex coordinate functions $z_1,\ldots, z_N$, and we denote by $R(X)$ the uniform closure of the rational functions  holomorphic on (a neighborhood of) $X$.  It is well known that the maximal ideal space of $P(X)$ can be naturally identified with the \emph{polynomially convex hull} $\what X$ of $X$ defined by
$$\what X=\{z\in\C^N:|p(z)|\leq \max_{x\in X}|p(x)|\
\mbox{\rm{for\ all\ polynomials}}\ p
\},$$
and the maximal ideal space of $R(X)$ can be naturally identified with the \emph{rationally convex hull} $\rhull$ of $X$ defined by
$$\rhull = \{z\in\C^N: p(z)\in p(X)\ 
\mbox{\rm{for\ all\ polynomials}}\ p
\}.$$

An equivalent formulation of the definition of $\rhull$ is that $\rhull$ consists precisely of those points $z\in \C^N$ such that every polynomial that vanishes at $z$ also has a zero on $X$.

For $E$ a subset of $\C^N=\R^{2N}$, by $m(E)$ we will denote the $2N$-dimensional Lebesgue measure of $E$.
The real part of a complex number (or function) $z$ will be denoted by $\Re z$.  In Sections~5 and~6 the following notation will be used in which we assume that $X$ is a compact set in $\C^N$ and that $\Omega$ is an open set in $\C^N$ each of which contains the origin:
\hfil\break
$\Bone X$ denotes the set of rational functions $f$ holomorphic on a neighborhood of $X$ such that $\|f\|_X\leq 1$,\hfil\break
$\Bzero X$ denotes the set of functions in $\Bone X$ that vanish at the origin,\hfil\break
$\Bone \Omega$ denotes the set of functions $f$ holomorphic on $\Omega$ such that \break
$\|f\|_\Omega\leq 1$, and\hfil\break
$\Bzero \Omega$ denotes the set of functions in $\Bone \Omega$ that vanish at the origin.\hfil\break

Let $A$ be a uniform algebra on a compact space $X$.
A point $x\in X$ is said to be a \emph{peak point} for $A$ if
there exists $f \in A$ with $f(x) = 1$ and $|f(y)| < 1$ for all $y \in X \setminus \{x\}$.
The \emph{Gleason parts} for the uniform algebra $A$ are the equivalence classes in the maximal ideal space of $A$ under the equivalence relation $\varphi\sim\psi$ if $\|\varphi-\psi\|<2$ in the norm on the dual space $A^*$.  (That this really is an equivalence relation is well-known but {\it not\/} obvious!)
We say that a Gleason part is \emph{nontrivial} if it contains more than one point.
For $\phi$ a multiplicative linear functional on $A$, a \emph{point derivation} on $A$ at $\phi$ is a linear functional $\psi$ on $A$ satisfying the identity 
$$\phantom{\hbox{for all\ } f,g\in A.} \psi(fg)=\psi(f)\phi(g) + \phi(f)\psi(g)\qquad \hbox{for all\ } f,g\in A.$$
A point derivation is said to be \emph{bounded} if it is bounded (continuous) as a linear functional.
It is well known that every peak point is a one-point Gleason part, and that at a peak point there are no nonzero point derivations.  (See for instance \cite{Browder}.)

An \emph{analytic disc} in the maximal ideal space $\cma$ of a uniform algebra $A$ is, by definition, an injective map $\sigma$ of the open unit disc $D\subset \C$  into 
$\cma$ such that the function $f\circ \sigma$ is analytic on $D$ for every $f$ in 
$A$.  It is immediate that the presence of an analytic disc implies the existence of a nontrivial Gleason part and nonzero bounded point derivations.

\font\tenscr=rsfs10 
\font\sevenscr=rsfs7 
\font\fivescr=rsfs5 
\skewchar\tenscr='177 \skewchar\sevenscr='177 \skewchar\fivescr='177
\newfam\scrfam \textfont\scrfam=\tenscr \scriptfont\scrfam=\sevenscr
\scriptscriptfont\scrfam=\fivescr
\def\scr{\fam\scrfam}

We will make repeated use of the following standard lemma.  (For a proof see \cite[Lemma~2.6.1]{Browder}.)

\begin{lemma}
Two multiplicative linear functionals $\phi$ and $\psi$ on a uniform algebra $A$ lie in the same Gleason part if and only if
$$\sup\{|\psi(f)|: f\in A, \|f\|\leq 1, \phi(f)=0\}< 1.$$
\end{lemma}

We will also need the following theorem of Wilken~\cite{Wilken} (or see \cite[Theorem~3.3.7]{Browder}).

\begin{theorem}
If $X$ is a compact set in the plane and $x\in X$ is not a peak point for $R(X)$, then the Gleason part of $x$ has positive planar measure.
\end{theorem}

It follows from Wilken's theorem that $R(X)$ ($X\subset \C$) has a Gleason part of positive planar measure whenever $R(X)\not= C(X)$, since in that case some point of $X$ must be a nonpeak point by a well-known theorem of Errett Bishop \cite{Bishop} (or see \cite[Theorem~3.3.3]{Browder}).


We will also need an elementary observation about bounded point derivations. Let $X\subset\C^n$ be compact.  Each bounded point derivation on $R(X)$ is of course determined by its restriction to $R_0(X)$, the space of rational functions holomorphic on a neighborhood of $X$.  It is not difficult to show that each point derivation on $R_0(X)$ is determined by its action on the functions $z_1,\ldots, z_N$, and hence, is a linear combination of $\partial /\partial z_1, \ldots, \partial/\partial z_N$.  We conclude that the complex vector space of bounded point derivation on $R(X)$ at a particular point $x\in X$ has dimension at most $N$, and that the dimension is exactly $N$ if and only if there is a number $M<\infty$ such that
\beaa
|\pfz(x)|\leq M \quad \hbox{for every } f\in \Bone {X} \hbox{\ and\ } l=\range N. \label{XDbound-n}
\eeaa

Finally we recall the standard method for translating 
statements about rational hulls in $\C^N$ into statements about polynomial hulls in $\C^{N+1}$.  A theorem, due to Kenneth Hoffman and Errett Bishop in case $N=1$, and Hugo Rossi~\cite{Rossi} in general, asserts that for $X$ a compact set in $\C^n$, there is a function $g$ such that the functions $z_1,\ldots, z_N, g$ generate $R(X)$ as a uniform algebra.  Furthermore, as noted by Basener~\cite{Basener}, the function $g$ can be taken to be $C^\infty$ on $\C^N$.  Let $\tau:\C^N\rightarrow \C^{N+1}$ be given by $g(z)=g\bigr(z, g(z)\bigl)$.  Then $P(\tau(X))$ is isomorphic as a uniform algebra to $R(X)$.  Consequently, $\widehat {\tau(X)} =g(\rhull)$.


\section{The theorem of Duval and Levenberg\\ and its analogue for rational convexity} 

In this section we prove Theorem~\ref{rationalDL}, the analogue for rational convexity of the result of Duval and Levenberg \cite{DuvalL} on polynomially convex hulls without analytic discs.  We also partially reprove the result of Duval and Levenberg.  We do this for several reasons: (i) to be able to explain clearly how these two proofs are related to each other and how they are both related to an earlier argument of Wermer, (ii) to make available a key lemma from the argument of Duval and Levenberg, as it is this lemma rather than the theorem itself that will be needed in the proof of Theorem~\ref{maintheorem1} , (iii) to establish that the set $X$ constructed by Duval and Levenberg satisfies $\wX=\rhull$ as claimed in the introduction (Theorem~\ref{DL}), and (iv) to give a simplification of a part of  the argument that is also needed in the proof Theorem~\ref{rationalDL}.

Throughout this section it is to be understood that that $B$ is the open unit ball in $\C^N$ with $N\geq 2$.

The proof of the following lemma is an easy exercise.

\blem \label{pconvex}
Let $X$ be a polynomially convex set in $\C^N$, let $Y$ be a polynomially convex set in the plane, and let $f$ be a polynomial on $\C^N$.  Then $f^{-1}(Y)\cap X$ is polynomially convex in $\C^N$.\hfill$\square$
\elem

The next lemma is really the foundation for the proof of Duval and Levenberg.

\blem \label{found}
Let $p$ be a polynomial on $\C^N$ and $X=\{\Re p \leq 0\}\cap \pB$.  Then $\wX=\rhull=\{\Re p\leq 0\}\cap\ob$.
\elem

\bpf
First note that for every real number $\alpha$, the sets $\{\Re p\leq \alpha\} \cap \ob$ and $\{\Re p\geq \alpha\} \cap \ob$ are polynomially convex by the preceding lemma.  Therefore,
$$\rhull \subset \wX\subset \{\Re p\leq 0\} \cap \ob.$$
To see that the reverse implications hold, suppose $z_0$ is a point of $\ob$ such that $z_0\notin \rhull$.   Then there is a polynomial $q$ such that $q(z_0)=0$ and $q$ has no zeros on $X$.  Thus $\Re p>0$ everywhere on the set $\{q=0\}\cap \pB$.  Hence there is an $\alpha>0$ such that $\{q=0\}\cap \pB \subset \{\Re p\geq \alpha\} \cap \pB$.  Now the maximum principle applied to the irreducible component of $\{q=0\}$ through $z_0$ gives that $z_0$ lies in the polynomial hull of $\{\Re p\geq\alpha\} \cap \pB$ and hence lies in $\{\Re p\geq\alpha\} \cap \ob$.  Thus $z_0\notin 
\{\Re p\leq 0\} \cap \ob$.
\epf

The key construction of Duval and Levenberg is contained in the next lemma.

\blem \label{keylemma}
If $K\subset B$ is a compact polynomially convex set, and if $\{p_j\}_{j=1}^\infty$ is a sequence of polynomials such that each of the sets $Z_j=\ob\cap \pji(0)$ is disjoint from $K$, then there is a compact subset $X$ of $\pB$ such that $\wX=\rhull\supset K$ and  $\wX\cap Z_j=\emptyset$ for every $j$.
\elem

\bpf
One constructs a sequence of polynomials $\{f_j\}_{j=1}^\infty$ such that for each $j$ we have $\Re f_j<0$ on $K$ and $\Re f_j>0$ on $Z_j$ and such that the sets 
\[X_j=\{\Re f_j\leq 0\}\cap \pB\]
form a decreasing sequence.  The set
$X=\bigcap_{j=1}^\infty X_j$
then has the required properties by the preceding lemma.
We refer the reader to the original paper \cite{DuvalL} or to \cite{Stout} for the construction of the $f_j$.
\epf

Finally we offer the following simplification of the remainder of the proof of Duval and Levenberg.

\bpf[Proof of Theorem~\ref{DL}]
The collection of all polynomials on $\C^N$ having no zeros on $K$, viewed as a subset of $P(\ob)$, has a countable dense subset.  Choosing such a subset $\{p_j\}$ and applying Lemma~\ref{keylemma} yields a compact subset $X$ of $\pB$ such that $\wX=\rhull\supset K$ and each $p_j$ is zero-free on $\wX$.  Assume to get a contradiction that $\wX\setminus K$ contains an analytic disc $\sigma:D\rightarrow \wX\setminus K$.  Since $K$ is polynomially convex, and hence rationally convex, there is a polynomial $p$ such that $p(\sigma(0))=0$ and $p$ has no zeros on $K$.  Because $\sigma$ is injective, we may assume, by adding to $p$ a small multiple of a suitable first degree polynomial if necessary,
that  $\bigl(\partial (p\circ \sigma)/\partial z\bigr)(0)\not= 0$.  Then Rouche's theorem shows that every polynomial that is uniformly close to $p$ on $\sigma(D)$ also has a zero on $\sigma(D)$.  But then some $p_j$ must have a zero on $\sigma(D)\subset \wX$, a contradiction. 
\epf

Duval and Levenberg refer to their construction as a modification of the construction of Stolzenberg.  However, the present author sees their construction as being more closely related to the construction of a {\em rationally} convex hull with no analytic discs due to Wermer~\cite{Wermerpark} (or see~\cite[Chapter~24]{AW}).
In both the construction of Duval and Levenberg and the construction of Wermer, the set having hull without analytic discs is obtained from the boundary of a domain by removing a sequence of subsets, and presence of nontrivial hull arises because a point lies in the rational hull of the set that remains if it does not lie in the polynomial hull of the set removed.  (Compare the proofs of Lemma~\ref{found} above and \cite[Claim~2, p.~207]{AW}.)  The principal difference between the proofs is that in Wermer's argument, the sets removed are inverse images of discs, while in Duval and Levenberg, they are inverse images of half-planes.  This difference is what accounts for the polynomial hull, and not only the rational hull, having no analytic discs in the construction of Duval and Levenberg.  Another difference between the constructions is that in the case of Duval and Levenberg, to get the example, one must consider general polynomials, whereas in the construction of Wermer, it suffices to remove inverse images of discs under only the functions $z_1$ and $z_2$.  These observations motivate our proof of Theorem~\ref{rationalDL}, the analogue for rational convexity of the Duval and Levenberg result.  In fact, our proof is, roughly, a generalization of Wermer's argument with $z_1$ and $z_2$ replaced by a sequence of polynomials.  In connection with this, note that the construction of Dales and Feinstein~\cite{DalesF} is a generalization of the construction of Stolzenberg also involving replacing $z_1$ and $z_2$ by a sequence of polynomials.


For the proof of Theorem~\ref{rationalDL} we need lemmas analogous to those above.

\blem \label{Wermerlemma}
Suppose $\{p_j\}$ is a sequence of polynomials and $\{r_j\}$ is a sequence of strictly positive numbers.
Define 
$$X=\pB\setminus \bigcup_{j=1}^\infty \{z\in \partial B: |p_j(z)|< r_j\}$$
 and for each $j$,
$$H_j=\{z\in \partial B: |p_j(z)|< r_j\}.$$
For $z_0\in B$, if $z_0\notin \rhull$, then for some $n$, we have $z_0\in [\cup_{j=1}^n \oh_j]\endhat$.
\elem

\bpf
We follow the argument of Wermer given in \cite[pp.~207--208]{AW}.  Suppose $z_0\notin \rhull$.  Then there is a polynomial $q$ such that $q(z_0)=0$ and $q$ has no zeros on $X$.  
Then $\{z\in \pB:q(z)=0\}\subset \cup_{j=1}^\infty H_j$.  since $\{z\in \pB: q(z)=0\}$ is compact, there is an $n$ such that $\{z\in \pB:q(z)=0\}\subset \cup_{j=1}^n H_j$.

Let $V$ be the irreducible component of the zero set of $q$ containing $z_0$.  Then $V\cap \pB\subset \cup_{j=1}^\infty H_j$. The maximum principle applied to $V$ gives that for every polynomial $p$,
$$|p(z_0)|\leq \max_{V\cap\pB} |p| \leq \max_{\bigcup_1^n \oh_j} |p|.$$
Hence $z_0\in [\cup_{j=1}^n \oh_j]\endhat$, as asserted.
\epf

\blem \label{keylemmarat}
If $K\subset B$ is a rationally convex set, and if $\{p_j\}$ is a sequence of polynomials such that all the sets $Z_j=\ob\cap \pji(0)$ are disjoint from $K$, then there is a compact subset $X$ of $\partial B$ such that $\rhull \supset K$ and $\rhull \cap Z_j=\emptyset$ for all $j$.
\elem

Note that the assertion that $\rhull \cap Z_j=\emptyset$ is equivalent to the statement that each $p_j$ is zero-free on $X$.

\bpf 
We will show that for a suitable choice of strictly positive numbers $r_1, r_2, \ldots$, the set $X$ defined as in Lemma~\ref{Wermerlemma} satisfies $\rhull\supset K$.  Since it is obvious that each $p_j$ is zero-free on $X$, this will prove the lemma.

We will choose the $r_j$ inductively so that the following statement holds:
for each $n=1, 2, \ldots$, there exists a polynomial $q_n$ so that, with $H_j$ defined as in Lemma~\ref{Wermerlemma}, 
\be
\sup_{z\in \bigcup_1^n \oh_j} |q_n(z)| <  \inf_{z\in K} |q_n(z)|. \label{inductioninequality}
\ee
Then $\rhull\supset K$ by Lemma~\ref{Wermerlemma}.

Since $p_1$ has no zeros on $K$, there is an $r_1>0$ such that $|p_1|>r_1$ on $K$.  Set $q_1=p_1$.  Then (\ref{inductioninequality}) holds with $n=1$.  Now assume that strictly positive numbers $\row rn$ and polynomials $\row qn$ have been found so that (\ref{inductioninequality}) holds.  Note that by multiplying $q_n$ by a suitable constant, we may assume that 
\beaa
\inf_{z\in K} |q_n(z)|=1.
\eeaa
Since then $\sup_{\cup_1^n \oh_j} |q_n|<1$, and $\inf_{K}|p_{n+1}|>0$, choosing $k$ large enough, we get
\beaa
\left(\sup_{\cup_1^n \oh_j} |q^k_n|\right) \left(\sup_{\cup_1^n \oh_j} |p_{n+1}|\right )<\inf_K|p_{n+1}|.
\eeaa
Choose $r_{n+1}>0$ small enough that, with $H_{n+1}$ defined as in Lemma~\ref{Wermerlemma},
\beaa
\left(\sup_{\oh_{n+1}} |q_n^k|\right) r_{n+1} < \inf_K |p_{n+1}|.
\eeaa
Now defined the polynomial $q_{n+1}$ by $q_{n+1}=q_n^k\cdot  p_{n+1}$.  Then (\ref{inductioninequality}) holds with $n$ replaced by $n+1$.  This completes the induction and the proof.
\epf

\bpf [Proof of Theorem~\ref{rationalDL}]
This follows from Lemma~\ref{keylemmarat} by exactly the same argument that was used above to obtain Theorem~\ref{DL} of Duval and Levenberg from Lemma~\ref{keylemma}.
\epf


\section{Hulls with dense invertibles}

In this section we prove Theorems~\ref{maintheorem1}--\ref{exp}.  We begin by noting the following general result (a modification of the result of Duval and Levenberg) giving hulls with dense invertibles.

\bthm \label{DLmod}
If $K\subset B\subset\C^N$ {\rm(}$N\geq 2${\rm)} is a compact polynomially convex set, and if the set of polynomials that are zero-free on $K$ is dense in $P(\ob)$,
then 
there is a compact subset $X$ of $\pB$ such that $\wX=\rhull\supset K$ and  $P(X)$ has dense invertibles.
\ethm

\bpf
The hypotheses enable us to choose a countable dense subset  $\{ p_j\}$ of $P(\ob)$ consisting of polynomials  each of which is zero-free on $K$.  The result then follows immediately from Lemma~\ref{keylemma}.
\epf

\blem \label{positive}
Let $\{p_j\}$ be a sequence of polynomials on $\C^N$, and let $\tB$ be a closed ball in $\C^N$.  There exists a compact polynomially convex set $E\subset \tB$ of positive $2N$-dimensional measure such that each $p_j$ is zero-free on $E$. 
\elem

\bpf
For each $j$, the set $\pji(\{z\in \C: \Re z=0\})$ is a real-analytic variety in $\C^N$ and hence has $2N$-dimensional measure zero.  Consequently,
$$m\Bigl(\pji(\{z\in \C: |\Re z|<\varep\})\cap \tB\Bigr) \rightarrow 0 \ \hbox{as\ }\varep \rightarrow 0.$$
Let $c=m(\tB)/2$.  For each $j$, choose $\varep_{j}>0$ such that
$$m\Bigl(\pji(\{z\in \C: |\Re z|<\varep_{j}\})\cap \tB\Bigr)< c/2^{j}.$$
Then
$$m\Bigl(\bigcup_{j=1}^\infty \pji(\{z\in \C: |\Re z|<\varep_{j}\})\cap \tB\Bigr)<\sum_{j=1}^\infty c/2^{j}=c.$$
Thus setting $E=\tB\setminus \bigcup_{j=1}^\infty \pji(\{z\in \C: |\Re z|<\varep_{j}\})$ we have $m(E)>c$.  Also that $E$ is compact and each $p_j$ is zero-free on $E$. 

For each $j$, choose a closed disc $\od_j$ containing $p_j(\tB)$.  Then
$$E=\bigcap_{j=1}^\infty \Bigl(\pji (\{z\in \od_j: |\Re z|\geq \varep_{j}\})\cap \tB \Bigr).$$
Each set $\{z\in \od_j: |\Re z|\geq \varep_{j}\}$ is  polynomially convex since it has connected complement in the plane.  Hence $\pji (\{z\in \od_j: |\Re z|\geq \varep_{j}\})\cap \tB$ is polynomially convex by Lemma~\ref{pconvex}.  Consequently, $E$ is  polynomially convex.
\epf

The proof of the next lemma is an easy application of the inverse function theorem together with Sard's theorem.

\blem \label{sard}
Let $f:\R^N\rightarrow \R^N$ be a smooth map, and let $E\subset \R^N$ be a closed set with empty interior.  Then $f(E)$ has empty interior in $\R^N$. \phantom{mmm}\hfill
$\square$
\elem

\bpf[Proof of Theorem~\ref{maintheorem1}]
Choose a Swiss cheese $L$ such that $R(L)$ has a bounded point derivation at almost every point.  (One can use the Swiss cheese of Browder given in \cite{Wermerderivation} or the Swiss cheese of Theorem~\ref{cheese} of the present paper.)  By rescaling, we may assume that $L$ is contained in the open unit disc.  As discussed in Section~2, there is a $C^\infty$ function $g$ on the plane such that the functions $z$ and $g$ generate $R(L)$ as a uniform algebra.  Let $\tau:\C\rightarrow \C^2$ be defined by $\tau(z)=\bigl(z,g(z)\bigr)$, and let $K=\tau(L)$.  By rescaling $g$, we may assume that $K$ is contained in $B$.  Note that $K$ is polynomially convex (see Section~2).    For every polynomial $p$ on $\C^2$, the map $p\circ \tau:\C\rightarrow\C$ is $C^\infty$, and hence $p(K)=(p\circ\tau)(L)$ has empty interior in the plane by Lemma~\ref{sard}.  Consequently, the collection of polynomials with no zeros on $K$ is dense in $P(\ob)$ and hence contains a countable subcollection $\{p_j\}$ that is also dense in $P(\ob)$.

Fix a nonconstant polynomial $q$.  
Because $q(K)$ is compact while $q(B)$ is open in $\C$, there is a closed disc $D_K$ that contains $q(K)$ but does not contain $q(B)$.  Consequently, there is a closed ball $\tB$ contained in $B$ such that $q(\tB)$ is contained in a closed disc $D_{\tB}$ disjoint from $D_K$.
By Lemma~\ref{positive}, contained in $\tB$ there is a compact polynomially convex set $E$ of positive $2N$-dimensional measure
such that each $p_j$ is zero-free on $E$.  Because $q(K)$ and $q(E)$ are contained in the disjoint discs $D_K$ and $D_{\tB}$, their polynomially convex hulls are disjoint.  Therefore, 
$K\cup E$ is polynomially convex by Kallin's lemma \cite{Kallin} (or see \cite[Theorem~1.6.19]{Stout}).  Also each $p_j$ is zero-free on $K\cup E$.  Theorem~\ref{DLmod} now yields the existence of a compact subset $X$ of $\pB$ such that $\wX=\rhull\supset K\cup E$ and  $P(X)$ has dense invertibles.

Recall (see Section~2) that $R(L)$ must have a Gleason part of positive 
2-dimensional measure, which we will denote by $Q_L$.  Furthermore, by our choice of $L$, there is a nonzero bounded point derivation for $R(L)$ at every point of some subset $P_L$ of full measure in $Q_L$.  Let $Q=\tau(Q_L)$, and let $P=\tau(P_L)$.  Then $Q$, and hence $P$, is contained in a single Gleason part for $P(X)$.  Furthermore,  $P$ has positive 2-dimensional Hausdorff measure (since $z_1(P)=P_L$).

Note that if $p$ is a polynomial on $\C^2$, then $p\circ \tau$ is in $R(L)$.  Consequently, if $f$ is a function in $P(X)$, and we denote the Gelfand transform of $f$ (regarded as a function on $\wX$) by $\hat f$, then $\hat f\circ \tau$ is in $R(L)$.  The map $P(X)\rightarrow R(L)$ given by $f\mapsto \hat f\circ \tau$ is clearly a norm-decreasing homomorphism with dense range.  It follows that if $\psi_L$ is a nonzero bounded point derivation on $R(L)$ at a point $z_0$, then the formula
$$\psi(f)=\psi_L(\hat f\circ \tau)$$
defines a nonzero bounded point derivation on $P(X)$ at $\tau(z_0)$.
Consequently, there is a nonzero bounded point derivation for $P(X)$ at every point of $P$.
\epf

\begin{remark}
As pointed out by Dales, using the Swiss cheese constructed by Anthony O'Farrell in \cite{OFarrellinf} in place of the one used above, we obtain a set $X\subset \C^N$ with nontrivial polynomial hull such that $P(X)$ has dense invertibles and a bounded derivation of infinite order.
\end{remark}

We now turn to rational hulls.  For the proof of Theorem~\ref{maintheorem3} we will need, in place of the Swiss cheese used in the proof of Theorem~\ref{maintheorem1}, a set constructed in Section~6 where Theorem~\ref{cheese-dim-n} is proved.

We begin with the analogue of Theorem~\ref{DLmod}.

\bthm \label{rationalDLmod}
If $K\subset B\subset \C^N$ {\rm(}$N\geq 2${\rm)} is a compact rationally convex set, and if the set of polynomials that are zero-free on $K$ is dense in $P(\ob)$,
then 
there is a compact subset $X$ of $\pB$ such that $\rhull\supset K$ and  $R(X)$ has dense invertibles.
\ethm

\bpf
The proof is similar to the proof of Theorem~\ref{DLmod}.  The hypotheses enable us to choose a countable dense subset  $\{ p_j\}$ of $P(\ob)$ consisting of polynomials  each of which is zero-free on $K$.  The result then follows immediately from Lemma~\ref{keylemmarat} and the following lemma.
\epf

\blem \label{poly-to-rat}
Let $X\subset \C^N$ be compact.  Suppose there is a dense subset of $P(X)$ consisting of polynomials each of which is zero-free on $X$.  Then $R(X)$ has dense invertibles.
\elem

\bpf
Consider an arbitrary rational function $r=p/q$ with $q$ zero-free on $X$.  Choose a sequence of polynomials $(p_j)$ with no zeros on $X$ such that $p_j\rightarrow p$ in $P(X)$.  Then $p_j/q\rightarrow p/q$ in $R(X)$.  Since each function $p_j/q$ is zero-free on $\rhull$, this proves the theorem.
\epf

\bpf [Proof of Theorem~\ref{maintheorem3}]
Let $K$ denote the set $X$ whose existence is given by 
Theorem~\ref{strongcheese}, and let $P$ be as in that theorem.  Then take $X$ to be the set given by Theorem~\ref{rationalDLmod}.  That conditions (i) and (ii) hold is clear.  To verify condition (iii), let $T:R(X) \rightarrow R(K)$ be the map sending each function in $R(X)$ to its restriction to $K$.  This map is clearly a 
norm-decreasing homomorphism with dense range.  Therefore, the adjoint map $T^*: R(K)^*\rightarrow R(X)^*$ is injective.  Furthermore, $T^*$ sends point derivations to point derivations.  Condition (iii) follows.
\epf

\bpf[Proof of Theorem~\ref{maintheorem2}]
This follows immediately from Theorem~\ref{maintheorem3} using the method discussed in the last paragraph of Section~2.
\epf

\bpf[Proof of Theorem~\ref{exp}]
The set $\inv A \setminus \exp A$ is open in any commutative Banach algebra.  (Here $\inv A$ denotes the set of invertible elements of $A$.) Thus it suffices to show that the algebras in question can be chosen so that $\exp A\subsetneq\inv A$.

Consider a circle $\{a+re^{i\theta}: 0\leq\theta\leq 2\pi\}$ contained in the Swiss cheese $L$ used in the proof of Theorem~\ref{maintheorem1}.  We can choose the collection of polynomials $\{p_j\}$ to contain the function $z_1-a$  Then $z_1-a$ is invertible in $P(X)$ but has no continuous logarithm on $X$.

To handle the algebra in Theorem~\ref{maintheorem3}, let $C$ be the circle $\{(\textfrac{1}{2} + \textfrac{1}{4}e^{i\theta},0,\ldots,0)\in \C^N: 0\leq\theta \leq 2\pi\}$.  Since the image of $C\cup \{0\}$ under any polynomial has empty interior in the plane, the set of polynomials with no zeros on $C\cup\{0\}$ is dense in $P(\ob)$.  Thus we can choose the collection $\{p_j\}$ in the proof of Theorem~\ref{strongcheese} to consist of polynomials that have no zeros on $C\cup\{0\}$ and to contain the function $z_1-\textfrac{1}{2}$.  Then the $r_j$ in the proof of Theorem~\ref{strongcheese} can be chosen so that $C$ is contained in the rationally convex set constructed there.  The set $C$ is then in the set $X$ obtained in the proof of Theorem~\ref{maintheorem3}, and the function $z_1-\textfrac{1}{2}$ is in $\inv R(X)\setminus \exp R(X)$.  

The assertion about the algebra in Theorem~\ref{maintheorem2} follows immediately from the case of Theorem~\ref{maintheorem3} since the two algebras in question are isomorphic.
\epf

\section{A Swiss cheese with a Gleason part of full measure} \label{cheesesection}

This section is devoted to the proof of Theorem~\ref{cheese}.  

Recall that we denote the open unit disc in the plane by $D$.  The open disc of radius $r$ with center $a$ will be denoted by $D(a,r)$, and the corresponding closed disc will be denoted by $\od(a,r)$.  The discs $D(0,r)$ and $\od(0,r)$ will be denoted by $rD$ and $r\od$ respectively.  Also recall the notations 
$\Bone X$,
$\Bzero X$,
$\Bone \Omega$, and
$\Bzero \Omega$ introduced in Section~2.

The proof of Theorem~\ref{cheese} will use several lemmas.

\blem\label{functionbound}
Suppose $K$ is a compact set containing the origin and contained in an open set $\Omega\subset \C^N$.  Then there exists an $R$ with $0<R<1$ such that $\|f\|_K\leq R$ for all $f\in \Bzero \Omega$.
\elem

\bpf
This is proven by a normal families argument left to the reader.
\epf

\blem\label{derivativebound}
Suppose $K$ is a compact set contained in an open set $\Omega\subset \C^N$.  Then there exists an $M<\infty$  such that $\|\pfz\|_K\leq M$ for all $f\in \Bone \Omega$ and $l=\range N$.
\elem

\bpf
This follows from Cauchy's estimates.
\epf

\blem\label{functioncontrol}
Let $K$ be a compact set containing the origin and contained in an open set $\Omega\subset \C$, let $a\in \Omega\setminus K$, and let $\varep>0$ be given.  Let $0<R<1$, and suppose that $\|f\|_K\leq R$ for all $f\in \Bzero \Omega$.
Then there exists an $r>0$ such that $\odisc ar{} \subset \Omega\setminus K$ and $\|f\|_K\leq R+\varep$ for all $f\in \Bzero {\Omega\setminus \odisc ar{}}$.
\elem

\bpf
Assume to get a contradiction that no such $r>0$ exists.  Let $J$ be a positive integer large enough that $\od(a, 1/J)\subset \Omega\setminus K$.  Then for each $n=J, J+1, \ldots$, there exists a function $f_n\in \Bzero {\Omega\setminus \od(a, 1/n)}$ such that $\|f_n\|_K> R+\varep$.  For each $m=J, J+1, \ldots$, the set $\{f_n: n\geq m\}$ is a normal family on $\Omega\setminus \od(a, 1/m)$.
Thus there is a subsequence of $(f_n)$ that converges uniformly on compact subsets of $\Omega \setminus \od(a, 1/J)$.  We may then choose a further subsequence that converges uniformly on compact subsets of  $\Omega \setminus \od\bigl(a, 1/(J+1)\bigr)$.  Continuing in this manner taking subsequences of subsequences, and then applying the usual diagonalization argument, we arrive at a subsequence $(f_{n_k})$ of $(f_n)$ such that for each compact subset $L$ of $\Omega\setminus \{a\}$, the sequence  $(f_{n_k})$
converges uniformly on $L$.  Note that given $L$, there may be a finite number of terms of the sequence $(f_{n_k})$ that are not defined on $L$, but this does not matter.  Thus there is a well-defined limit function $f$ that is holomorphic on $\Omega\setminus \{a\}$.  Since $\|f\|_{\Omega\setminus \{a\}}\leq 1$, the singularity at $a$ is removable, so we may regard $f$ as defined on $\Omega$.  Then $f$ is in $\Bzero\Omega$.  Thus by hypothesis, $\|f\|_K\leq R$.  But  since $f_{n_k}\rightarrow f$ uniformly on $K$, and $\|f_n\|_K>R+\varep$ for all $n$, this is a contradiction.
\epf

\blem\label{derivativecontrol}
Let $K$ be a compact set contained in an open set $\Omega\subset \C$, let $a\in \Omega\setminus K$, and let $\varep>0$ be given.  Let $0<M<\infty$, and suppose that $\|f'\|_K\leq M$ for all $f\in \Bone \Omega$.
Then there exists an $r>0$ such that $\odisc ar{} \subset \Omega\setminus K$ and $\|f'\|_K\leq M+\varep$ for all $f\in \Bone {\Omega\setminus \odisc ar{}}$.
\elem

\bpf
This is proven by a normal families argument similar to the one just presented and hence is left to the reader.
\epf

\blem\label{measurezero}
Let $\{a_j\}$ be a sequence in $D$, and let $\{r_j\}$ and $\{u_j\}$ be sequences of real numbers such that $0<r_j<u_j$ for every $j$, and $\sum u_j<\infty$.  Set 
$$X=\od\setminus \bigcup\limits_{j=1}^\infty D(a_j,r_j)$$
and
$$K_n=\left(1-{\textfrac{1}{n+1}}\right)\od\setminus \left(\left[\bigcup_{\,j=1}^{n-1} \disc a{(1+{\textfrac{1}{n}})r}j \right] \cup \left[\bigcup_{\,j=n}^\infty \disc auj\right]\right).$$
Then $\bigcup_{n=1}^\infty K_n\subset X$ and $m(X\setminus \bigcup_{n=1}^\infty K_n)=0$.
\elem

\bpf
The asserted inclusion is obvious.
Now note that for each $k=1,2, \ldots$,
\begin{eqnarray*}
X\setminus \bigcup_{n=1}^\infty K_n &\subset& X\setminus K_k \\ 
&\subset &
\left(\od\setminus \left(1-{\textfrac{1}{k+1}}\right)\od\right)
\cup \left(\bigcup_{j=1}^{k-1} \Bigl[ \disc a{(1+{\textfrac{1}{k}})r}j)\setminus \disc arj\Bigr]\right)\\
& & \qquad\cup \left(\bigcup_{j=k}^{\infty} \disc auj\right).
\end{eqnarray*}
Thus
$$m(X\setminus \bigcup_{n=1}^\infty K_n) \leq \pi(1-(1-{\textfrac{1}{k+1}})^2) +\sum_{j=1}^{k-1} \pi ((1+{\textfrac{1}{k}})^2 -1) r_j^2 + \sum_{j=k}^\infty \pi u_j^2.$$
Since the right hand side goes to zero as $k\rightarrow \infty$, the left hand side must be zero, as desired.
\epf

It might be tempting to imagine that $X\setminus \bigcup_{n=1}^\infty K_n = \bigcup_{n=1}^\infty \partial D(a_j,r_j)$.  However, this is not in general the case, as $X\setminus \bigcup_{n=1}^\infty K_n$ is a $G_\delta$ set while $\bigcup_{n=1}^\infty \partial D(a_j,r_j)$ need not be.

\bpf[Proof of Theorem~\ref{cheese}]
We will choose a sequence $\{a_j\}$ in $D$ and sequences $\{r_j\}$ and $\{u_j\}$ of real numbers such that $0<r_j<u_j$ for every $j$, and $\sum u_j<\infty$ such that the set $X=\od \setminus \cup D(a_j,r_j)$ is a   Swiss cheese that contains the origin and for each set $K_n$ defined as in Lemma~\ref{measurezero}, there are constants $0<C_n<1$ and $0<C'_n<\infty$ such that for every $z\in K_n$ we have
\be
\sup_{f\in \Bzero {X}}|f(z)|\leq C_n\label{Xbound}
\ee
 and
\be
\sup_{f\in \Bone {X}}|f'(z)|\leq C'_n. \label{XDbound}
\ee
Then by facts about Gleason parts and bounded point derivations discussed in Section~2, (\ref{Xbound}) gives that every point of each $K_n$ lies in the Gleason part of the origin for $R(X)$, and (\ref{XDbound}) gives that there is a bounded point derivation on $R(X)$ at each point of each $K_n$.  Theorem~\ref{cheese} then follows by Lemma~\ref{measurezero} in view of the theorem of Wilken mentioned in Section~2.

First choose a sequence of points $\{\alpha_j\}$ dense in $D$.  Set $a_1=\alpha_1$, and choose $0<u_1<1/2$ small enough that $\odisc au1\subset D$ and $0\notin \odisc au1$.
Set 
$$K_1^1=(1/2) \od \setminus \disc au1.$$
By Lemmas~\ref{functionbound} and~\ref{derivativebound}
there exist $0<R_1<1$ and $0<M_1<\infty$ such that 
\[
\|f\|_{K_1^1}\leq R_1 \quad \forallBzero D
\]
and
\[
\|f'\|_{K_1^1}\leq M_1 \quad \forallBone D
\]
Now Lemmas~\ref{functioncontrol} and~\ref{derivativecontrol} give that there exists $0<r_1<u_1$ such that 
\[
\|f\|_{K_1^1}\leq R_1 + \textfrac 14 (1-R_1) \quad \forallBzero {D\setminus \odisc ar1}
\]
and
\[
\|f'\|_{K_1^1}\leq (1+\textfrac12) M_1 \quad \forallBone {D\setminus \odisc ar1}.
\]
Of course we can choose $r_1$ so that  $\partial \disc ar1$ is disjoint from the set $\{\alpha_j\}$

Next set $a_2=\alpha_v$ where $v$ is the smallest positive integer such that $\alpha_v\notin \disc ar1$.  Choose $0<u_2<1/2^2$ small enough that $\odisc au2\subset D\setminus\odisc ar1$ and $0\notin \odisc au2$. Set
$$K_1^2=K_1^1\setminus \disc au2=(1/2) \od \setminus [\disc au1 \cup \disc au2]$$
and
$$K_2^2=(2/3) \od \setminus [\disc a{{\textfrac 32} r}1 \cup \disc au2].$$ 
Note that each of $K_1^2$ and $K_2^2$ is contained 
in $D\setminus \odisc ar1$ and neither contains the point $a_2$.
Lemmas~\ref{functionbound} and~\ref{derivativebound} yield numbers $0<R_2<1$ and $0<M_2<\infty$ such that
\[
\|f\|_{K_2^2}\leq R_2 \quad \forallBzero {D\setminus \odisc ar1}
\]
and
\[
\|f'\|_{K_2^2}\leq M_2 \quad \forallBone {D\setminus \odisc ar1}.
\]
Lemmas~\ref{functioncontrol} and~\ref{derivativecontrol} now yield the existence of $0<r_2<u_2$ such that
\[
\|f\|_{K_1^2}\leq R_1 + (\textfrac 14 + \textfrac 18) (1-R_1) \quad \forallBzero {D\setminus [ \odisc ar1 \cup \odisc ar2]}
\]

\[
\|f\|_{K_2^2}\leq R_2 + \textfrac 14 (1-R_2) \quad \forallBzero {D\setminus [ \odisc ar1 \cup \odisc ar2]}
\]

\[
\|f'\|_{K_1^2}\leq (1+\textfrac12 +\textfrac 14) M_1 \quad \forallBone {D\setminus [ \odisc ar1 \cup \odisc ar2]}
\]
and
\[
\|f'\|_{K_2^2}\leq (1+\textfrac12) M_2 \quad \forallBone {D\setminus [ \odisc ar1 \cup \odisc ar2]}.
\]
Of course we can choose $r_2$ so  that $\partial \disc ar2$ is disjoint from the set $\{\alpha_j\}$.

We then continue by induction.  Suppose that for some $k\geq 2$, we have chosen $a_1,\ldots, a_k$, $u_1,\ldots, u_k$, and $r_1,\ldots, r_k$, and for $1\leq n\leq m\leq k$ we have set
\be
K_n^m=\left(1-{\textfrac{1}{n+1}}\right)\od\setminus \left(\left[\bigcup_{\,j=1}^{n-1} \disc a{(1+{\textfrac{1}{n}})r}j \right] \cup 
\left[ \bigcup_{\,j=n}^m \disc auj\right]\right), \label{Kdef}
\ee
and for all $1\leq l\leq k$ and all $1\leq n\leq m\leq k$
the following conditions hold:
\begin{enumerate}
\item[(i)]  $a_l$ is the first $\alpha_v$ such that $\alpha_v\notin \disc ar1\cup\cdots\cup \disc ar{l-1}$
\item[(ii)] $0<r_l<u_l<1/2^l$
\item[(iii)] $0\notin \odisc aul\subset D\setminus [\odisc ar1 \cup\cdots\cup \odisc ar{l-1}]$
\item[(iv)] $\partial \disc arl$ is disjoint from the set $\{\alpha_j\}$
\item[(v)] there exists $0<R_n<1$ such that 
$$\|f\|_{K_n^m}\leq R_n +(\textfrac 14 + \textfrac 18 + \cdots + \textfrac 1{2^{m-n+2}})(1-R_n) \quad \forallBzero {D\setminus \cup_{j=1}^m \odisc arj}$$
\item[(vi)] there exists $M_n<\infty$ such that 
$$\|f'\|_{K_n^m}\leq (1+ \textfrac 12 +\cdots + \textfrac 1{2^{m-n+1}})M_n \quad \forallBone {D\setminus \cup_{j=1}^m \odisc arj}.$$
\end{enumerate}
Choose $a_{k+1}$ to be the first $\alpha_v$ such that $\alpha_v\notin \disc ar1\cup\cdots\cup \disc ar{k}$.  Choose $0<u_{k+1}<1/2^{k+1}$ small enough that $0\notin \odisc au{k+1}\subset D\setminus \cup_{j=1}^k \odisc arj$.  For $n=1,\ldots, k+1$, define $K_{n}^{k+1}$ so that (\ref{Kdef}) continues to hold with $m$ replaced by $k+1$.  Note that then for $n=1,\ldots, k$, we have 
\[
K_n^{k+1}=K_n^k\setminus \disc au{k+1}.
\]
Note also that each of $K_1^{k+1}, \ldots, K_{k+1}^{k+1}$ is contained 
in $D\setminus \cup_{j=1}^k \odisc arj$ and none contains the point $a_{k+1}$.

Lemmas~\ref{functionbound} and~\ref{derivativebound} yield numbers $0<R_{k+1}<1$ and $0<M_{k+1}<\infty$ such that
\[
\|f\|_{K_{k+1}^{k+1}}\leq R_{k+1} \quad \forallBzero {D\setminus \cup_{j=1}^k \odisc arj}
\]
and
\[
\|f'\|_{K_{k+1}^{k+1}}\leq M_{k+1} \quad \forallBone {D\setminus \cup_{j=1}^k \odisc arj}.
\]
Lemmas~\ref{functioncontrol} and~\ref{derivativecontrol}, together with conditions (v) and (vi) of the induction hypothesis, now yield the existence of $0<r_{k+1}<u_{k+1}$ such that
conditions (v) and (vi) continue to hold when, in the restriction $1\leq n\leq m\leq k$,  we replace $k$ by $k+1$. 
Of course we can choose $r_{k+1}$ so  that $\partial \disc ar{k+1}$ is disjoint from the set $\{\alpha_j\}$.  Thus the induction can proceed.

We conclude that we can obtain sequences $\{a_j\}$, $\{r_j\}$, and $\{u_j\}$ such that with $K_n^m$ defined by (\ref{Kdef}) for all $1\leq n\leq m$, we have for all $1\leq l$ and all $1\leq n\leq m$ that conditions (i)--(vi) hold.  Note that the sets $K_n$ defined as in Lemma~\ref{measurezero} satisfy
$K_n=\bigcap_{m=n}^\infty K_n^m$.

Finally consider an arbitrary point $z\in K_n$.  Given $g\in \Bzero X$, there is some $m\geq n$ such that $g\in \Bzero {D\setminus \cup_{j=1}^m \odisc arj}$.  Since $z\in K_n^m$, condition (v) above gives that 
$$|g(z)|\leq R_n +(\textfrac 14 + \textfrac 18 + \cdots + \textfrac 1{2^{m-n+2}})(1-R_n).$$  Therefore,
$$\sup_{f\in \Bzero X} |f(z)| \leq R_n + \textfrac 12 (1-R_n) <1$$
so (\ref{Xbound}) holds with $C_n=R_n + \textfrac 12 (1-R_n)$.
Similarly, using condition (vi) we get that
$$\sup_{f\in \Bone X} |f'(z)| \leq 2M_n<\infty$$
so (\ref{XDbound}) holds with $C'_n=2M_n$.
This completes the proof.
\epf

\begin{remark}
As mentioned in the introduction, $R(X)$ is antisymmetric for the Swiss cheese $X$ just constructed.  This is immediate from the easy proposition below because the Gleason part of the origin is dense in $X$.  The denseness of of the Gleason part of the origin is easily seen from the construction above.  However,  it is in fact easily shown that for any Swiss cheese, every subset of full measure is dense, as we now demonstrate.  Given a Swiss cheese $K$, let $\mu$ be the measure that is $dz$ on the boundary of the unit disc and $-dz$ on the boundary of each removed disc.  Given an arbitrary open set $U$ of the plane, let $h$ be a nonnegative smooth function on $\C$ with compact support contained in $U$.  Then, letting $\hat\mu$ denote the Cauchy transform of $\mu$, the measure $h\, d\mu - (1/\pi)(\partial h /\partial {\overline z}) \hat\mu\, dm$ is a nonzero measure with support in $U$ that annihilates $R(K\cap \overline U)$.  (See \cite[Section II.10]{Gamelin}.)The Hartogs-Rosenthal theorem thus implies that no (relatively) open subset of $K$ can have planar measure zero. 
\end{remark}

\begin{proposition}
Suppose $A$ is a uniform algebra on a compact space $X$ and $f$ is a 
real-valued function in $A$.  Then $f$ is constant on each Gleason part for $A$.
\end{proposition}

\bpf
Suppose $x$ and $y$ are points such that $f(x)\not=f(y)$.  Let $g:\R\rightarrow [-1,1]$ be a continuous function such that $g\bigl(f(x)\bigr)=1$ and $g\bigl(f(y)\bigr)=-1$.  By the Weierstrass approximation theorem $g$ can be approximated uniformly on $f(X)$ by polynomials.  Hence the function $g\circ f$ is in $A$.  Since $(g\circ f)(x)=1$, $(g\circ f)(y)=-1$, and $\|g\circ f\|=1$, we conclude that $x$ and $y$ lie in different Gleason parts.
\epf


\section{Proof of Theorem~\ref{cheese-dim-n}} \label{cheese-dim-n-section}
This section is devoted to the proof of Theorem~\ref{cheese-dim-n}.

Recall that $D$ denoted the open unit disc in the plane, and $rD$ denotes the disc of radius $r$ centered at the origin.  Recall also the notations
$\Bone X$,
$\Bzero X$,
$\Bone \Omega$, and
$\Bzero \Omega$ introduced in the introduction.

We will obtain Theorem~\ref{cheese-dim-n} by proving the following ostensibly stronger result.  The condition that the polynomials that are zero-free on $X$ are dense in $P(\ob)$ was used earlier in the proof of Theorem~\ref{maintheorem3}

\bthm \label{strongcheese}
There exists a compact rationally convex set $X\subset B\subset\C^N$ {\rm (}$N\geq 1${\rm )} of positive $2N$-dimensional measure  such that  the collection of polynomials zero-free on $X$ is dense in $P(\ob)$  and such that there is a set $P\subset X$ of full $2N$-dimensional  measure in $X$ such that $P$ is contained in a single Gleason part for $R(X)$ and at every point of $P$ the space of bounded point derivations for 
$R(X)$ has dimension $n$.

\ethm

\bpf [Proof of Theorem~\ref{cheese-dim-n}]
This is immediate from Theorem~\ref{strongcheese} and Lemma~\ref{poly-to-rat}.
\epf

The basic idea and overall structure of the proof of Theorem~\ref{strongcheese} is the same as the proof of Theorem~\ref{cheese} although there are differences in the details.  First we give modifications of some lemmas from Section~5.

\blem\label{functioncontrol-rat}
Let $K$ be a compact set containing the origin and contained in an open set $\Omega\subset \C^N$, let $p$ be a polynomial on $\C^N$ with no zeros on $K$, and let $\varep>0$ be given.  Let $0<R<1$, and suppose that $\|f\|_K\leq R$ for all $f\in \Bzero \Omega$.
Then there exists an $r>0$ such that $p^{-1}(r\od)\cap K=\emptyset$ and $\|f\|_K\leq R+\varep$ for all $f\in \Bzero {\Omega\setminus p^{-1}(r\od)}$.
\elem

\bpf
This is proved by a normal families argument similar to the proof of the next lemma and hence is left to the reader. 
\epf

\blem\label{derivativecontrol-rat}
Let $K$ be a compact set contained in an open set $\Omega\subset \C^N$, let $p$ be a polynomial on $\C^N$ with no zeros on $K$, and let $\varep>0$ be given.  Let $0<M<\infty$, and suppose that $\| \pfz \|_K\leq M$ for all $f\in \Bone \Omega$ and $l=\range N$.
Then there exists an $r>0$ such that $p^{-1}(r\od) \cap K=\emptyset$ and $\|\pfz\|_K\leq M+\varep$ for all $f\in \Bone {\Omega\setminus p^{-1}(r\od)}$ and $l=\range N$.
\elem

\bpf
Of course it suffices to show for each fixed $l=\range N$, there exists an $r_l>0$ such that $p^{-1}(r_l\od) \cap K=\emptyset$ and $\|\pfz\|_K\leq M+\varep$ for all $f\in \Bone {\Omega\setminus p^{-1}(r_l\od)}$, so let $l$ be fixed and assume to get a contradiction that no such $r_l$ exists.  

Let $J$ be a positive integer large enough that $p^{-1}((1/J)\od) \cap K=\emptyset$.  Then for each $n=J, J+1, \ldots$, there exists a function $f_n\in \Bzero {\Omega\setminus p^{-1}((1/J)\od)}$ such that $\|\partial f_n / \partial z_l\|_K> M+\varep$.  For each $m=J, J+1, \ldots$, the set $\{f_n: n\geq m\}$ is a normal family on $\Omega\setminus p^{-1}((1/m)\od)$.
Thus arguing as in the proof of Theorem~\ref{functioncontrol} taking subsequences of subsequence and applying the usual diagonalization, we can obtain
a subsequence $(f_{n_k})$ of $(f_n)$ such that for each compact subset $L$ of $\Omega\setminus p^{-1}(0)$, the sequence  $(f_{n_k})$
converges uniformly on $L$.  Note that given $L$, there may be a finite number of terms of the sequence $(f_{n_k})$ that are not defined on $L$, but this does not matter.  Thus there is a well-defined limit function $f$ that is holomorphic on $\Omega\setminus p^{-1}(0)$.  Since $\|f\|_{\Omega\setminus p^{-1}(0)}\leq 1$, 
the extended Riemann removable singularities theorem \cite[Theorem~D2]{Gunning} gives that $f$ has a holomorphic extension to $\Omega$, which we continue to denote by $f$.
Then $f$ is in $\Bone\Omega$.  Thus by hypothesis, $\|\pfz\|_K\leq M$.  But  since $\partial f_{n_k} / \partial z_l\rightarrow \partial f / \partial z_l$ uniformly on $K$, and $\|\partial f_n / \partial z_l\|_K>M+\varep$ for all $n$, this is a contradiction.
\epf

\blem \label{measurezero-dim-n}
Let $\{p_j\}$ be a sequence of nonconstant polynomials on $\C^N$, and let $\{r_j\}$ and $\{u_j\}$ be sequences of real numbers such that $0<r_j<u_j/2$ for every $j$ and such that $\sum_{j=1}^\infty m\Bigl(\pji(u_jD)\cap \ob\Bigr)<\infty$.
Set 
$$X=\ob\setminus \bigcup\limits_{j=1}^\infty \pji(r_jD)$$
and
$$K_n=\left(1-{\textfrac{1}{n+1}}\right)\ob\setminus \left(\left[\bigcup_{\,j=1}^{n-1} 
\pji \bigl((1+{\textfrac{1}{n}})r_jD\bigr) \right] \cup 
\left[ \bigcup_{\,j=n}^\infty \pji(u_jD)\right]\right).$$
Then $\bigcup_{n=1}^\infty K_n\subset X$ and $m(X\setminus \bigcup_{n=1}^\infty K_n)=0$.
\elem

\bpf
The asserted inclusion is obvious.  Now it suffices to show that $m(X\setminus K_n)\rightarrow 0$ as $n\rightarrow \infty$.  Note that
\begin{eqnarray*}
X\setminus K_n &\subset&
\bigl(\ob \setminus (1-\textfrac{1}{n+1}) \ob\bigr)\cup
\left(\ob\cap \left[\bigcup_{\,j=1}^{n-1} \pji \Bigl( (1+\textfrac{1}{n})r_j D \setminus r_jD\Bigr)\right]\right) \\
&&\qquad\cup \left(\ob \cap \left[\bigcup_{\,j=n}^\infty \pji (u_j D)\right]\right).
\end{eqnarray*}
Let 
$$c_{n,j} =  \begin{cases} m\Bigl(\ob \cap \pji \bigl( (1+\textfrac{1}{n})r_j D \setminus r_jD\bigr)\Bigr)     &{\rm for \ \ } j\leq n-1 \cr  m\Bigl(\ob \cap \pji (u_jD)\Bigr)  & {\rm for \ \ }  j\geq n. \end{cases} $$
Then 
\begin{eqnarray} \label{volumebound}
m(X\setminus K_n) \leq \bigl( 1- (1-\textfrac{1}{n+1})^{2N}\bigr) m(\ob) + \sum_{j=1}^\infty c_{n,j}.
\end{eqnarray}
Note that
$$ \ob \cap \left[ \bigcap_{n=1}^\infty \pji \bigl( (1+\textfrac{1}{n})r_j D \setminus r_jD\bigr)\right] =\ob \cap \pji \bigl(\partial (r_jD) \bigr)
$$
and
$$m\bigl(\pji \bigl(\partial (r_jD) \bigr)=0$$
since $\pji \bigl(\partial (r_jD)$ is a real-analytic subvariety of $\C^N$ of real dimension strictly less than $2N$.
Consequently, for each fixed $j$, we have $c_{n,j}\rightarrow 0$ as $n\rightarrow\infty$.  Thus we conclude (by the dominated convergence theorem for instance) that  $\sum_{j=1}^\infty c_{n,j} \rightarrow 0$ as $n\rightarrow \infty$.  Thus inequality (\ref{volumebound}) implies that $m(X\setminus K_n)\rightarrow 0$, as desired.
\epf

With the lemmas in hand, we proceed to the proof of Theorem~\ref{strongcheese}.

\bpf[Theorem~\ref{strongcheese}]
Choose a countable collection $\{p_j\}$ of nonconstant polynomials on $\C^N$ such that $\|p_j\|_B\leq 1$ and $p_j(0)\not =0$ for every $j$ and such that the set $\{\lambda p_j: \lambda \in \C, j=1,2,\ldots\}$ is dense in $P(\ob)$.  We will choose strictly positive numbers $r_1, r_2, \ldots$ such that the set $X$ defined by
$$X=\ob\setminus \bigcup\limits_{j=1}^\infty \pji(r_jD)$$
has the desired properties.  Fix $\varep>0$.  We will in fact arrange for $m(X)\geq m(\ob)-\varep$.

Some of the properties are automatic.  It is obvious that $X$ is a compact rationally convex set.  Since each $p_j$ is zero-free on $X$, the collection of polynomials zero-free on $X$ is dense in $P(\ob)$.

We will choose sequences $\{r_j\}$ and $\{u_j\}$ of positive numbers such that the hypotheses of Lemma~\ref{measurezero-dim-n} are satisfied and so that  for each set $K_n$ defined as in Lemma~\ref{measurezero-dim-n}, there are constants $0<C_n<1$ and $0<C'_n<\infty$ such that for every $z\in K_n$ we have
\be
\sup_{f\in \Bzero {X}}|f(z)|\leq C_n, \label{Xbound-n}
\ee
 and
\be
\sup_{f\in \Bone {X}}|(\pfz)(z)|\leq C'_n \quad \hbox{for\ all\ } l=\range N. \label{XDbound-n}
\ee
Then by facts about Gleason parts and bounded point derivations discussed in Section~2, (\ref{Xbound-n}) gives that every point of each $K_n$ lies in the Gleason part of the origin for $R(X)$, and (\ref{XDbound-n}) gives that at each point of each $K_n$ the space of bounded point derivations on $R(X)$ has dimension $N$.  Theorem~\ref{strongcheese} then follows by Lemma~\ref{measurezero-dim-n}.

Note that each set $\pji(0)$ is an analytic subvariety of $\C^N$ with no interior and hence has measure zero.  Consequently, $m\bigl(\pji(rD)\bigr)\rightarrow 0$ as $r\rightarrow 0$.  Also $0\notin \pji(rD)$ for $r$ sufficiently small.

Choose $u_1>0$ small enough that $m\bigl(\pjid {u_1}\bigr)<\varep/2$ and $0\notin \pjid {u_1}$.
Set 
$$K_1^1=(1/2) \ob \setminus p_1^{-1}(u_1D).$$
By Lemmas~\ref{functionbound} and~\ref{derivativebound}
there exist $0<R_1<1$ and $0<M_1<\infty$ such that 
\[
\|f\|_{K_1^1}\leq R_1 \quad \forallBzero B
\]
and
\[
\|\pfz\|_{K_1^1}\leq M_1 \quad \forallBone B \hbox{\ and\ } l=\range N.
\]
Now Lemmas~\ref{functioncontrol-rat} and~\ref{derivativecontrol-rat} give that there exists $0<r_1<u_1/2$ such that 
\[
\|f\|_{K_1^1}\leq R_1 + \textfrac 14 (1-R_1) \quad \forallBzero {B\setminus p_1^{-1}(r_1\od)}
\]
and
\[
\|\pfz \|_{K_1^1}\leq (1+\textfrac12) M_1 \quad \forallBzero {B\setminus p_1^{-1}(r_1\od)} \hbox{\ and\ } l=\range N.
\]

Next choose $u_2>0$ small enough that $m\bigl(p_2^{-1}(u_2D)\bigr)<\varep/2^2$ and $0\notin p_2^{-1}(u_2D)$. Set
$$K_1^2=K_1^1\setminus p_2^{-1}(u_2D)=(1/2) \ob \setminus [p_1^{-1}(u_1D)\cup p_2^{-1}(u_2D)]$$
and
$$K_2^2=(2/3) \ob \setminus [p_1^{-1}({\textfrac 32} r_1 D) \cup p_2^{-1}(u_2 D)].$$ 
Note that each of $K_1^2$ and $K_2^2$ is contained 
in $B\setminus p_1^{-1}(r_1 D)$.
Lemmas~\ref{functionbound} and~\ref{derivativebound} yield numbers $0<R_2<1$ and $0<M_2<\infty$ such that
\[
\|f\|_{K_2^2}\leq R_2 \quad \forallBzero{B\setminus p_1^{-1}(r_1D)}
\]
and
\[
\|\pfz\|_{K_2^2}\leq M_2 \quad \forallBone {B\setminus p_1^{-1}(r_1D)} \hbox{\ and\ } l=\range N.
\]
Lemmas~\ref{functioncontrol-rat} and~\ref{derivativecontrol-rat} now yield the existence of $0<r_2<u_2/2$ such that
\[
\|f\|_{K_1^2}\leq R_1 + (\textfrac 14 + \textfrac 18) (1-R_1) \quad \forallBzero 
{B\setminus [ p_1^{-1}(r_1D) \cup p_2^{-1}(r_2D)]}
\]

\[
\|f\|_{K_2^2}\leq R_2 + \textfrac 14 (1-R_2) \quad \forallBzero 
{B\setminus [ p_1^{-1}(r_1D) \cup p_2^{-1}(r_2D)]}
\]

\[
\|\pfz\|_{K_1^2}\leq (1+\textfrac12 +\textfrac 14) M_1 \quad \forallBone 
{B\setminus [ p_1^{-1}(r_1D) \cup p_2^{-1}(r_2D)]} \hbox{\ and\ } l=\range N
\]
and
\[
\|\pfz\|_{K_2^2}\leq (1+\textfrac12) M_2 \quad \forallBone 
{B\setminus [ p_1^{-1}(r_1D) \cup p_2^{-1}(r_2D)]} \hbox{\ and\ } l=\range N.
\]

We then continue by induction.  Suppose that for some $k\geq 2$, we have chosen $u_1,\ldots, u_k$ and $r_1,\ldots, r_k$, and for $1\leq n\leq m\leq k$ we have set
\be \label{Kdef-rat}
K_n^m=\left(1-{\textfrac{1}{n+1}}\right)\ob\setminus \left(\left[\bigcup_{\,j=1}^{n-1} 
\pji \bigl((1+{\textfrac{1}{n}})r_jD\bigr) \right] \cup 
\left[ \bigcup_{\,j=n}^m \pji(u_jD)\right]\right)
\ee
and for all $1\leq l\leq k$ and all $1\leq n\leq m\leq k$
the following conditions hold:
\begin{enumerate}
\item[(i)] $0<r_l<u_l/2$
\item[(ii)] $m\bigl(p_l^{-1}(u_lD)\bigr)<\varep/2^l$
\item[(iii)] $0\notin p_l^{-1}(u_lD)$
\item[(iv)] there exists $0<R_n<1$ such that 
$$\|f\|_{K_n^m}\leq R_n +(\textfrac 14 + \textfrac 18 + \cdots + \textfrac 1{2^{m-n+2}})(1-R_n) \quad \forallBzero {B\setminus \cup_{j=1}^m p_1^{-1}(r_1\od)}$$
\item[(v)] there exists $M_n<\infty$ such that 
$$\|\pfz\|_{K_n^m}\leq (1+ \textfrac 12 +\cdots + \textfrac 1{2^{m-n+1}})M_n \quad \forallBone {B\setminus \cup_{j=1}^m p_1^{-1}(r_1\od)} \hbox{\ and\ } l=\range N.$$
\end{enumerate}
Next choose $u_k>0$ small enough that $m\bigl(p_k^{-1}(u_kD)\bigr)<\varep/2^k$ and $0\notin p_k^{-1}(u_kD)$. For $n=1,\ldots, k+1$, define $K_{n}^{k+1}$ so that (\ref{Kdef-rat}) continues to hold with $m$ replaced by $k+1$.  Note that then for $n=1,\ldots, k$, we have 
\[
K_n^{k+1}=K_n^k\setminus p_{k+1}^{-1}(u_{k+1}D).
\]
Note also that each of $K_1^{k+1}, \ldots, K_{k+1}^{k+1}$ is contained 
in $D\setminus \bigcup_{j=1}^k p_j^{-1}(r_j \od)$.

Lemmas~\ref{functionbound} and~\ref{derivativebound} yield numbers $0<R_{k+1}<1$ and $0<M_{k+1}<\infty$ such that
\[
\|f\|_{K_{k+1}^{k+1}}\leq R_{k+1} \quad \forallBzero {B\setminus\cup_{j=1}^k p_j^{-1}(r_j \od)}
\]
and
\[
\|\pfz\|_{K_{k+1}^{k+1}}\leq M_{k+1} \quad \forallBone {B\setminus\cup_{j=1}^k p_j^{-1}(r_j \od)} \hbox{\ and\ } l=\range N.
\]
Lemmas~\ref{functioncontrol-rat} and~\ref{derivativecontrol-rat}, together with conditions (v) and (vi) of the induction hypothesis, now yield the existence of $0<r_{k+1}<u_{k+1}/2$ such that
conditions (iv) and (v) continue to hold when, in the restriction $1\leq n\leq m\leq k$,  we replace $k$ by $k+1$. 
Thus the induction can proceed.

We conclude that we can obtain sequences $\{r_j\}$ and $\{u_j\}$ such that with $K_n^m$ defined by (\ref{Kdef-rat}) for all $1\leq n\leq m$, we have for all $1\leq l$ and all $1\leq n\leq m$ that conditions (i)--(v) hold.  Note that the sets $K_n$ defined as in Lemma~\ref{measurezero-dim-n} satisfy
$K_n=\bigcap_{m=n}^\infty K_n^m$.

Finally consider an arbitrary point $z\in K_n$.  Given $g\in \Bzero X$, there is some $m\geq n$ such that $g\in \Bzero {D\setminus \cup_{j=1}^m p_j^{-1}(r_j\od)}$.  Since $z\in K_n^m$, condition (iv) above gives that 
$$|g(z)|\leq R_n +(\textfrac 14 + \textfrac 18 + \cdots + \textfrac 1{2^{m-n+2}})(1-R_n).$$  Therefore,
$$\sup_{f\in \Bzero X} |f(z)| \leq R_n + \textfrac 12 (1-R_n) <1$$
so (\ref{Xbound-n}) holds with $C_n=R_n + \textfrac 12 (1-R_n)$.
Similarly, using condition (v) we get that
$$\sup_{f\in \Bone X} |(\pfz)(z)| \leq 2M_n<\infty \hbox{\ for\ all\ } l=\range N$$
so (\ref{XDbound-n}) holds with $C'_n=2M_n$.
This completes the proof.
\epf


\section{Open Questions}

\begin{enumerate}
\item[1.]  In Theorem~\ref{maintheorem1}, can $P$ be made to have positive 4-dimensional measure?  Can the space of bounded point derivations be made to have dimension 2 at each point of $P$?  In the conclusion of Theorem~\ref{maintheorem2}, can $N-1$ be replaced by $N$?
\item[2.]  Suppose that $X\subset \C^N$ is a compact set such that $\wX\setminus X$ is nonempty but contains no analytic discs.  Can it happen that $P(X)$ has no nonzero bounded point deriviations, but the Gleason part of some point of $\wX\setminus X$ is nontrivial?  Can it happen that $P(X)$ has no nontrivial Gleason parts, but there is a nonzero bounded point derivation at some point of $\wX \setminus X$?
\item[3.]  In connection with Wermer's Swiss cheese with no nonzero bounded point derivations and Theorem~\ref{cheese} we make the following conjecture: 
There exists a Swiss cheese $X$ such that the set of nonpeak points for $R(X)$ is a single Gleason part of full 2-dimensional measure in $X$, and $R(X)$ has no nonzero bounded point derivations.
\item[4.] Does there exist a Swiss cheese with a nonzero bounded point derivation at every nonpeak point?  One might be tempted to speculate that if two points lie in the same Gleason part and there is a nonzero bounded point derivation at one point then there must also be a nonzero bounded point derivation at the other point.  It would then follow at once that the Swiss cheese constructed above has the requested property.  However, such a speculation would be incorrect, as demonstrated by an example of O'Farrell~\cite{OFarrelliso} of a Swiss cheese with a nonzero bounded point derivation at exactly one point.  
\end{enumerate}

\end{document}